\documentclass[leqno]{amsart}

\usepackage{amsmath,amsfonts,amsthm,amssymb, mathrsfs, mathabx, accents, mathdots, textcomp, fontenc,  vntex}

\usepackage[french,english]{babel}

\usepackage[bookmarks]{hyperref}
\hypersetup{backref, colorlinks=true}

 \usepackage{graphicx}
\usepackage[all]{xy}
\newtheorem{theorem}{Theorem}[section]
\newtheorem{lemma}[theorem]{Lemma}

\newtheorem{remark}[theorem]{Remark}

\newtheorem{example}[theorem]{Example}
\newtheorem{definition}[theorem]{Definition}

\newtheorem{algorithm}[theorem]{Algorithm}

\newenvironment{equationth}{\stepcounter{theorem}\begin{equation}}{\end{equation}}
\newenvironment{preuve}{{\em{\noindent \textbf{Proof.} }}}
{\hfill $\blacksquare$}

\def\C{ \mathbb{C}}

\def\N{ \mathbb{N}}

\def\rond{\mathaccent"7017}

\begin{document}

\large

\title[]{ AN ALGORITHM TO CLASSIFY THE ASYMPTOTIC SET  \\ 
ASSOCIATED TO A POLYNOMIAL MAPPING}

\makeatother

\author[Nguy\~{\^e}n Th\d{i} B\'ich Th\h{u}y]{Nguy\~{\^e}n Th\d{i} B\'ich Th\h{u}y}
\address[{Nguy\~{\^e}n Th\d{i} B\'ich Th\h{u}y}]{UNESP, Universidade Estadual Paulista, ``J\'ulio de Mesquita Filho'', S\~ao Jos\'e do Rio Preto, Brasil}
\email{bichthuy@ibilce.unesp.br}
\maketitle \thispagestyle{empty}
\begin{abstract}
We provide an algorithm to classify the asymptotic sets of the dominant polynomial mappings $F: \C^3 \to \C^3$ of degree 2, using the definition of the so-called ``{\it fa{\c c}ons}'' in \cite{Thuy}. We obtain a  classification theorem for  the asymptotic sets of dominant polynomial mappings $F: \C^3 \to \C^3$ of degree 2. This algorithm can be generalized for the dominant polynomial mappings  $F: \C^n \to \C^n$ of degree $d$, with any $(n, d) \in {(\N^*)}^2$. 
\end{abstract}

\section{Introduction}
Let $F : \C^n \to \C^n$ be a polynomial mapping. Let us denote by $S_F$ the set of points at which  $F$ is non  proper, {\it i.e.}, 
$$S_F = \{ a \in \C^n \text{ such that } \exists \{ \xi_k\}_{k \in \N} \subset \C^n, \vert \xi_k \vert  \text{ tends to infinity and } F(\xi_k) \text{ tends to } a\},$$
where $ \vert \xi_k \vert$ is the Euclidean norm of  $\xi_k$ in $\C^n$. 
 The set $S_F$ is called the asymptotic set of $F$. In the years 90's, Jelonek studied this set in a deep way 
 and described the principal properties. One of the important results is that, if $F$ is dominant, {\it i.e.}, $\overline{F(\C^n)} = \C^n$, then $S_F$ is an empty set or a hypersurface \cite{Jelonek1}.

Notice that it is sufficient to define $S_F$  by considering sequences $\{\xi_k\}$ tending to infinity in the following sense: 
each coordinate of these sequences either tends to infinity or converges. In \cite{Thuy}, the sequences tending to infinity such that their  images tend to the points in $S_F$ are labeled in terms of  ``{\it fa{\c c}ons}'', as follows: 
for each point $a$ of $S_F$, there exists a sequence $\{\xi_k^a\}$ in the source space $\C^n$,  $\, \xi_k^a = (x_{k,1}^a, \ldots , x_{k,n}^a)$ 
tending to infinity such that 
 $F(\xi_k^a)$ tends to $a$. 
 Then there exists at least one index $i \in \N$, $1 \leq i \leq n$ 
such that  $x_{k,i}^a$ tends to infinity. 
We define a ``{\it fa{\c c}on''} of the point $a \in S_F$ 
 as a $(p,q)$-tuple $(i_1, \ldots , i_p)[j_1, \ldots, j_q]$ of integers 
where   
$x_{k,i_r}^{a}$ tends to infinity for $r = 1, \ldots , p$ and, 
for $s = 1, \ldots, q$, the sequence 
$ x_{k,j_s}^a$ tends to a complex number independently on the point $a$ 
 when $a$ describes locally $S_F$ (definition \ref{definitionXi}).  

The aim of this paper is to provide an algorithm to classify the asymptotic sets of dominant polynomial mappings $F: \C^3 \to \C^3$ of degree 2, using the definition of ``{\it fa{\c c}ons}'' in \cite{Thuy}, and then generalize this algorithm for the general case. 
 One important tool of the algorithm is  the notion of {\it pertinent variables}.  The idea of the notion of   pertinent variables is the following:  Let $F=(F_1, F_2 , F_3): \C^3  \to \C^3$ be  a dominant polynomial mapping  of degree 2  such that $S_F \neq \emptyset$. We fix a {\it fa{\c c}on} $\kappa$ of $F$ and assume that  $ \{ \xi_k \} $ is a sequence tending to infinity with the {\it fa{\c c}on} $\kappa$ such that $F(\xi_k)$ tends to a point of $S_F$. Since the degree of $F$ is 2 then each coordinate $F_1, F_2$ and $F_3$ of $F$ is a linear combination of $f_1 = x_1$, $f_2 = x_2$, $f_3 = x_3$, $f_4 = x_1x_2$, $f_5 = x_2x_3$ and $f_6 = x_3x_1$. We call a {\it pertinent variable of $F$ with respect to the {\it fa{\c c}on} $\kappa$}  a {\it minimum} linear combination   of $f_1, \ldots, f_6$ such that the image of the sequence $\{ \xi_k \}$ by this combination does not tend to infinity (see definition \ref{pertinentvar}). 


Moreover, if  $F$ is dominant then by Jelonek, the set $S_F$   has pure dimension $2$ (see theorem \ref{theoremjelonek1}). With this observation and with the idea of pertinent variables, we:

\begin{enumerate}
\item[$\bullet$] Make the list $(\mathcal{L})$  of all possible {\it fa{\c c}ons} 
 for a polynomial mappings $F: \C^3 \to \C^3$. 
This list is finite. In fact, there are 19 possible {\it fa{\c c}ons} (see the list (\ref{group})). 
\item[$\bullet$] Assume that a 2-dimensional irreductible stratum $S$ of $S_F$ admits $l$ fixed {\it fa{\c c}ons} in the list $(\mathcal{L})$, where $ 1 \leq l \leq 19$.
\item[$\bullet$] Determine the pertinent variables of $F$ with respect to these $l$ {\it fa{\c c}ons}. 
\item[$\bullet$] Restrict the above pertinent variables using the dominancy of $F$ and the fact $\dim S = 2$. We get the form of $F$ in terms of these pertinent variables.
\item[$\bullet$] Determine the geometry of $S$ in terms of the form of $F$.
\item[$\bullet$] Let $l$ runs in the list $(\mathcal{L})$ for $ 1 \leq l \leq 19$. 
We get all the possible $2$-dimensional irreductible strata of $S_F$. Since the dimension of $S_F$ is $2$, then we get the list of all possible 
asymptotic sets $S_F$. 
\end{enumerate}

With this idea,  we provide the algorithm \ref{algorithmordre} to classify the asymptotic sets
 of dominant polynomial mappings $F: \C^3 \to \C^3$ of degree 2, and 
we obtain the classification theorem  \ref{theothuyb}. 
This algorithm can be generalized for the general case of polynomial mappings $F: \C^n \to \C^n$ of degree $d$, where $n \geq 3$ and $d \geq 2$ (algorithm \ref{algorithmordregeneral}).

\section{Dominancy, assymptotic set and ``fa{\c c}ons''}
\subsection {Dominant polynomial mapping}
\begin{definition} \label{defdominant}
{\rm 
Let $F: \C^n \to \C^n$ be a polynomial mapping.  Let $\overline{F(\C^n)}$ be  the closure of $F(\C^n)$ in $\C^n$. F is called {\it dominant} if $\overline{F(\C^n)} = \C^n$, {\it i.e.},  $F(\C^n)$ is dense in $\C^n$. 
}
\end{definition}

We provide here a lemma on the dominancy of a polynomial mapping $F: \C^n \to \C^n$ that we will use later on. 


\begin{lemma}  \label{lemmaindependant}
Let $F =(F_1,\ldots, F_n) : \C^n \to \C^n$ be a dominant polynomial mapping. Then, the coordinate polynomials $F_1, \ldots, F_n$ are independent. That means, there does not exist any coordinate polynomial $F_\eta$, where $\eta \in \{1, \ldots, n \}$,  such that $F_\eta$ is a polynomial mapping of the variables $F_1, \ldots, F_{\eta-1}, F_{\eta+1}, \ldots, F_n$.
\end{lemma}
\begin{preuve}
 Assume that $F_\eta=\varphi (F_1, \ldots, F_{\eta-1}, F_{\eta+1}, \ldots, F_n)$ where  $\eta \in \{1, \ldots, n \}$ and
  $\varphi$ is a polynomial. Then,  the dimension of $F(\C^n)$ is  less than $n$. Consequently, the dimension of  $\overline{F(\C^n)}$ is less than $n$. That provides the contradiction with the fact $F$ is dominant. 
\end{preuve}

\subsection{Asymptotic set} \label{ensembleJelonek}
\begin{definition} \label{ensembleJelonek}
{\rm Let $F: \C^n \to \C^n$ be a polynomial mapping. Let us denote by $S_F$ the set of points at which  $F$ is non-proper, {\it i.e.}, 
\begin{equation*} \label{defSF}
S_F = \{ a \in \C^n \text{ such that } \exists \{ \xi_k\}_{k \in \N} \subset \C^n, \vert \xi_k \vert  \to \infty \text{  and } F(\xi_k) \to a\},
\end{equation*}
where $ \vert \xi_k \vert$ is the Euclidean norm of  $\xi_k$ in $\C^n$. 
The set $S_F$ is called the asymptotic set of $F$. 
}
\end{definition}

Recall that, it is sufficient to define $S_F$  by considering sequences $\{\xi_k\}$ tending to infinity in the following sense: 
each coordinate of these sequences either tends to infinity or converges to a finite number.

\begin{theorem} \cite{Jelonek1} \label{theoremjelonek1}
Let $F: \C^n \rightarrow \C^n$ be a polynomial mapping. 
 If $F$ is dominant, then $S_F$ is either an empty set or a hypersurface. 
\end{theorem}


\subsection{``Fa{\c c}ons''.} 
 
In this section, let us recall the definition of {\it fa{\c c}ons} as it appears in \cite{Thuy}. 
 In order to a better understanding of the definition of  {\it fa{\c c}ons}, let us start by giving an example. 
\begin{example} \cite{Thuy} \label{exfacon}
{\rm 
Let $F =  (F_1, F_2, F_3): \C^3_{(x_1, x_2, x_3)} \to \C^3_{(\alpha_1, \alpha_2, \alpha_3)}$ be the  polynomial mapping such that 
$$F_1:=x_1, \qquad F_2:= x_2, \qquad F_3:=x_1x_2x_3.$$
Notice that by the notations $\C^3_{(x_1, x_2, x_3)}$ and $ \C^3_{(\alpha_1, \alpha_2, \alpha_3)}$, we want to distinguish the source space and the target space. That means, if we take $x = (x_1, x_2, x_3)$ then $x$ belongs to the source space  $\C^3_{(x_1, x_2, x_3)}$; if we take $\alpha = (\alpha_1, \alpha_2, \alpha_3)$ then $\alpha$ belongs to the target space $ \C^3_{(\alpha_1, \alpha_2, \alpha_3)}$. 
We determine now the asymptotic set $S_F$ by using the definition \ref{ensembleJelonek}. 
 Assume that there exists a sequence 
$\{ \xi_k = ( x_{1, k}, x_{2, k}, x_{3,k}) \}$ in the source space  $\C^3_{(x_1, x_2, x_3)}$ 
tending to infinity such that its image $\{ F(\xi_k) = ( x_{1, k}, x_{2, k}, x_{1, k} x_{2, k} x_{3,k}) \}$ does not tend to infinity. Then $x_{1, k}$ and $x_{2, k}$ cannot tend to infinity. 
Since the sequence $\{\xi_k\}$ tends to infinity, then $x_{3, k}$ must tend to infinity. Hence, we have the three following cases: 

1)  $x_{1, k}$ tends to 0, $x_{2, k}$ tends to a complex number $\alpha_2 \in \C$ and $x_{3, k}$ 
 tends to infinity. In order to determine the biggest possible subset of $S_F$, 
we choose the sequences $x_{1, k}$ tending to 0 and $x_{3, k}$ tending to infinity in such a way that  the product 
$x_{1, k}x_{3, k}$ tends to a  complex number $\alpha_3$. 
 Let us choose, for example 
$\xi_k = \left( \frac{1}{k}, \alpha_2, \frac{k\alpha_3}{\alpha_2} \right)$ 
where $\alpha_2 \neq 0$, then $F(\xi_k)$ tends to a point $a = (0, \alpha_2, \alpha_3)$ in $S_F$. 
 We get a 2-dimensional stratum $S_1$ of $S_F$, where $S_1 = (\alpha_1 = 0) \setminus 0\alpha_3 \subset \C^3_{(\alpha_1, \alpha_2, \alpha_3)}$. 
We say that a {\it ``fa{\c c}on''} 
of $S_1$ is $(3)[1]$. 
The symbol ``(3)'' in the {\it fa{\c c}on}  $(3)[1]$ means that the {\it third} coordinate $ x_{3, k}$ of the sequence $\{ \xi_k \}$ tends to infinity. 
The symbol ``[1]'' in the {\it fa{\c c}on} $(3)[1]$ means that the {\it first} coordinate $x_{1, k}$ of the sequence $\{ \xi_k \}$  tends to 0 which is a fixed complex number which does not depend on the point $a = (0, \alpha_2, \alpha_3)$ when $a$ describes  $S_1$.  Notice that the second coordinate of the sequence $\{\xi_k\}$ tends to a complex number $\alpha_2$ depending on the point $a = (0, \alpha_2, \alpha_3)$ when $a$ varies, 
then the indice ``2'' does not appear in the {\it fa{\c c}on} $(3)[1]$.  Moreover, all the sequences tending to infinity such that their images tend to a point of $S_1$ admit only the {\it fa{\c c}on}  $(3)[1]$. 

The two following cases are  similar to the case 1): 

2) $x_{1, k}$ tends to a  complex number $\alpha_1 \in \C$, $x_{2, k}$ tends to 0 and $x_{3, k}$ 
 tends to infinity: then the {\it fa{\c c}on} $(3)[2]$ determines a 2-dimensional stratum $S_2$ of $S_F$, where $S_2 = (\alpha_2 = 0) \setminus 0\alpha_3 \subset \C^3_{(\alpha_1, \alpha_2, \alpha_3)}$. 

3) $x_{1, k}$ and $x_{2, k}$ tend to 0, and  $x_{3, k}$ 
 tends to infinity: then the {\it fa{\c c}on} $(3)[1, 2]$ determines the 1-dimensional stratum $S_3$ where $S_3$ is the axis $0\alpha_3$ in $\C^3_{(\alpha_1, \alpha_2, \alpha_3)}$. 

In conclusion, we get 
\begin{enumerate}
\item[$\bullet$]  
the asymptotic set $S_F$ of the given polynomial mapping $F$ as the union of two planes $(\alpha_1 = 0)$ and $(\alpha_2 = 0)$ in $\C^3_{(\alpha_1, \alpha_2, \alpha_3)}$, 
\item[$\bullet$] all the {\it fa{\c c}ons} of $S_F$ of the given polynomial mapping $F$: they are three fa{\c c}ons $(3)[1]$, $(3)[2]$ and $(3)[1, 2]$.
\end{enumerate}

\begin{remark}
{\rm The chosen sequence $\left\{ \xi_k = \left( \frac{1}{k}, \alpha_2, \frac{k\alpha_3}{\alpha_2} \right) \right\}$  in 1) of the above example is called a {\it generic sequence} of the 2-dimensional irreductible component $(\alpha_1 = 0)$ (a plane) of $S_F$, since the image of any sequence of this type (with differents $\alpha_2 \neq 0$ and $\alpha_3$) falls to a generic point of the plane $(\alpha_1 = 0)$. That means the images of all the sequences $\{\xi_k\}$ when $\alpha_2$ runs in $\C \setminus \{ 0 \}$ and $\alpha_3$ runs in $\C$  cover $S_1 = (\alpha_1 = 0) \setminus 0\alpha_3$ and $S_1$ is dense in the plane $(\alpha_1 = 0)$. We can see easily that a generic sequence of the 2-dimensional irreductible component $(\alpha_2 = 0)$ of $S_F$ is $\left( \alpha_1, \frac{1}{k}, \frac{k\alpha_3}{\alpha_1} \right)$ where $\alpha_1 \neq 0$.  
More generally, any sequence  $\left\{ \left( \frac{1}{k^r}, \alpha_2, \frac{k^s\alpha_3}{\alpha_2} \right) \right \}$, where $r = s$ and $\alpha_2 \neq 0$,  is a generic sequence of $(\alpha_1 = 0) \subset S_F$. 
Any sequence $\left \{ \left( \alpha_1, \frac{1}{k^r}, \frac{k^s\alpha_3}{\alpha_1} \right) \right \}$, where $r = s$ and $\alpha_1 \neq 0$, is a generic sequence of $(\alpha_2 = 0) \subset S_F$. 
} 
\end{remark}

}
\end{example}

In the light of this example, we recall here the definition of {\it fa{\c c}ons }  in \cite{Thuy}.

\begin{definition} \cite{Thuy} \label{definitionXi}
{\rm Let $F: \C^n \to \C^n$ be a dominant polynomial mapping  such that 
$S_F \neq \emptyset$.  
 For each point $a$ of $S_F$, there exists a sequence $\{ \xi_k^a \} \subset \C^n$,  $\, \xi_k^a = (x_{k,1}^a, \ldots , x_{k,n}^a)$ 
tending to infinity such that 
 $F(\xi_k^a)$ tends to $a$. 
 Then, there exists at least one index $i \in \N$, $1 \leq i \leq n$   
such that  $x_{k,i}^a$ tends to infinity when $k$ tends to infinity. 
 We define ``{\it a fa{\c c}on of tending to infinity of the sequence  $\{\xi_k^a\}$''},
  as a maximum $(p,q)$-tuple $\kappa = (i_1, \ldots , i_p)[j_1, \ldots, j_q]$ of different integers in $\{1, \ldots, n\}$, such that:
\begin{enumerate}
\item[i)] $x_{k,i_r}^{a}$  tends to infinity for all  $r = 1, \ldots , p$, 

\item[ii)] for all $s = 1, \ldots , q$, the sequence $x_{k,j_s}^a$  tends to a complex  number 
independently on the point $a$ 
 when $a$ varies locally, that means:
\begin{enumerate}
\item[ii.1)] either 
there exists in $S_F$ a subvariety
$U_a$ containing $a$ such that 
for any point $a'$ in $U_a$, 
there exists a sequence 
$\{ \xi_k^{a'} \} \subset \C^n$,  $\, \xi_k^{a'} =  (x_{k,1}^{a'}, \ldots , x_{k,n}^{a'})$ 
tending to infinity such that 
\begin{enumerate}
\item[a)] $F(\xi_k^{a'}) $ tends to $a',$
\item[b)] 
$x_{k,i_r}^{a'}$  tends to infinity for all  $r = 1, \ldots , p$, 
\item[c)] for all $s = 1, \ldots , q$,  
${\displaystyle \lim_{k \to \infty}} x_{k, j_s}^{a'} = {\displaystyle \lim_{k \to \infty}} x_{k,j_s}^{a}$ and this limit is finite.
\end{enumerate}

\item[ii.2)] or there does not exist such a subvariety, then we define 
$$\kappa = (i_1, \ldots , i_p)[j_1, \ldots, j_{n-p}],$$ 
where  $x_{k,i_r}^{a}$ tends to infinity for all $r = 1, \ldots , p$ and $\{i_1, \ldots , i_p\} \cup \{j_1, \ldots, j_{n-p}\} = \{1, \ldots, n\}$. 
In this case, the set of points $a$ is a subvariety of dimension 0 of $S_F$. 
\end{enumerate}
\end{enumerate}
We call  a {\it fa{\c c}on} of tending to infinity of the sequence  $\{\xi_k^a\}$ also a 
 {\it a fa{\c c}on}  of $S_F$. If the image of a sequence corresponding with a {\it fa{\c c}on}  $\kappa$ tends to a point of a stratum $S$ of $S_F$, we call also $\kappa$ a {\it fa{\c c}on}  of the stratum $S$. 

} 
\end{definition}


\section{An algorithm to stratify the asymptotic sets of the dominant polynominal mappings $F: \C^3 \to \C^3$ of degree 2}
 In this section we provide an algorithm to stratify the asymptotic sets associated to   dominant polynominal mappings $F: \C^3 \to \C^3$ of degree 2. In the last section, we show that this algorithm can be  generalized in the general case for  dominant polynominal mappings $F: \C^n \to \C^n$ of degree $d$ where $n \geq 3$ and $d \geq 2$. Let us recall that by degree of  a polynomial mapping  $F = (F_1, \ldots, F_n): \C^n \to \C^n$, we mean the highest degree of the monomials $F_1, \ldots ,  F_n$.

Let us consider now  a dominant polynomial mapping $F=(F_1, F_2 , F_3): \C^3_{(x_1, x_2, x_3)} \to \C^3_{(\alpha_1, \alpha_2, \alpha_3)}$ of degree 2  such that $S_F \neq \emptyset$. 
An important step of this section is to  define the notion of   {\it ``pertinent'' variables} of $F$. 
\subsection{Pertinent variables}

Let us explain at first the idea of the notion of {\it pertinent variables} of $F$: let $ \{ \xi_k \} = \{ (x_{1, k}, x_{2, k}, x_{3, k}) \}$ be a sequence in the source space $\C^3_{(x_1, x_2, x_3)}$ tending to infinity such that $F(\xi_k)$ tends to a point of $S_F$ in the target space $\C^3_{(\alpha_1, \alpha_2, \alpha_3)}$. Then the image of  $\xi_k$ by any coordinate polynomial  $F_\eta$, where $\eta = 1, 2, 3$, cannot tend to infinity. Notice that $F_\eta$ can be  written as  the  sum of elements of the form 
$F_\eta^1 F_\eta^2$ such that  
if 
 $F_\eta^1(\xi_k)$ tends to infinity, then $F_\eta^2(\xi_k)$ must tend to 0. 
In other words, if one element of the above sum  has a factor tending to  infinity with respect to the sequence $\{\xi_k\}$, then this element must be ``balanced'' with another factor tending to zero with respect to  the sequence $\{\xi_k\}$. 
 For example, assume that the coordinate sequences $x_{1, k}$ and $x_{2, k}$ of the sequence  $ \{ \xi_k \}$  tend to infinity, then  $F_\eta$  cannot admit neither $x_1$ nor $x_2$ {\it alone}  as an element of the above sum, but $F_\eta$ can admit 
$( x_1 - \nu x_2)$, $(x_1 - \nu x_2)  x_1$, $(x_1 - \nu x_2)  x_2$ 
 as elements of this sum, where $\nu \in \C \setminus \{ 0 \}$.  So we define


\begin{definition} \label{pertinentvar}
{\rm Let  $F=(F_1, F_2 , F_3): \C^3_{(x_1, x_2, x_3)} \to \C^3_{(\alpha_1, \alpha_2, \alpha_3)}$ be a polynomial mapping of degree 2 such that 
$S_F \neq \emptyset$. 
 Let us fix a {\it fa{\c c}on} $\kappa$ of $S_F$. 
Then there exists a sequence $\{ \xi_k \} \subset \C^3_{(x_1, x_2, x_3)}$ 
tending to infinity with the
 {\it fa{\c c}on} $\kappa$ such that its image tend to a point in $S_F$. An element   in the list 
\begin{equationth} \label{vairablepertinent}
\begin{cases}
 X_{h_i} = x_i, \text{ where } i = 1, 2, 3, \cr 
X_{h_j} = x_i + \nu_{h_j} x_j, \text{ where } i \neq j  \text{ and } i, j = 1, 2, 3, \cr 
 X_{h_r} = (x_i + \nu_{h_r} x_j) x_l, \text{ where } i \neq j \text{ and } i, j, l = 1, 2, 3, \cr 
 X_{h_s} = x_i +   \nu_{h_s} x_j x_l, \text{ where } i \neq j, j \neq k, k \neq i \text{ and } i, j, k = 1, 2, 3,
\end{cases}
\end{equationth} 
($ \nu_{h_i},  \nu_{h_j},  \nu_{h_r},  \nu_{h_s}  \in \C \setminus \{0\}$) 
 is called  
{\it a pertinent variable of $F$ with respect to the fa{\c c}on $\kappa$} 
if the image of the sequence $\{\xi_k\}$ by this element does not tend to infinity. 
}
\end{definition}

\begin{remark}
{\rm From now on, we will denote $X_1, \ldots, X_h$ pertinent variables of $F$  with respect to a {\it fa{\c c}on} and we write 
$$F = \tilde{F}(X_1, \ldots, X_h).$$

Notice that   we can also determine the pertinent variables of $F$ with respect to a set of  {\it fa{\c c}ons} in the case  we have more than one {\it fa{\c c}on}. }
\end{remark}

\subsection{Idea of the algorithm}

The aim of the algorithm that we present in this section is  to describe the list of all possible asymptotic sets $S_F$ for the  dominant polynomial mappings $F: \C^3 \to \C^3$ of degree 2. 
In order to do that, we observe firstly that 

\begin{enumerate}
\item[$\bullet$] The list  of all the possible {\it fa{\c c}ons} of $S_F$ for a polynomial mapping $F: \C^3 \to \C^3$ is 
\begin{equationth} \label{group}
\begin{cases}
\text{ 1) Group I}:  (1, 2, 3), \cr 
\text{ 2) Group II}: (1,2), (2,3) \text{ and } (3,1), \cr
\text{ 3) Group III}:  (1), (2) \text{ and }(3), \cr
\text{ 4) Group IV}: (1,2)[3], (1,3)[2] \text{ and } (2,3)[1], \cr 
\text{ 5) Group V}: (1)[2], (1)[3], (2)[1], [2](3), [3](1) \text{ and } [3](2), \cr
\text{ 6) Group VI}: (1)[2,3], (2)[1,3] \text{ and } (3)[1,2].
\end{cases}
\end{equationth}
This list have 19 {\it fa{\c c}ons}. 

\item [$\bullet$] Since $F$ dominant, then by the theorem \ref{theoremjelonek1}, the set $S_F$   has  pure dimension $2$.
\end{enumerate}

With these observations, we will: 

\begin{enumerate}
\item[$\bullet$] assume that a 2-dimensional irreductible stratum $S$ of $S_F$ admits $l$ fixed {\it fa{\c c}ons} in the list (\ref{group}), where $ 1 \leq l \leq 19$, 
\item[$\bullet$] determine the pertinent variables of $F$ with respect to these $l$ {\it fa{\c c}ons}, 
\item[$\bullet$] restrict the above pertinent variables using the dominancy of $F$ and the fact $\dim S = 2$. We get the form of $F$ in terms of these pertinent variables,
\item[$\bullet$] determine the geometry of $S$ in terms of  the form of $F$,
\item[$\bullet$] let $l$ runs in the list (\ref{group}) for $ 1 \leq l \leq 19$. 
We get all the possible $2$-dimensional irreductible strata $S$ of $S_F$. Since the dimension of $S_F$ is $2$, then we get the list of all the possible 
asymptotic sets $S_F$ of $F$. 
\end{enumerate}

The following example explains the process of the algorithm, {\it i.e. } how we can determine the geometry of a $2$-dimensional irreductible stratum $S$ of $S_F$ admitting some fixed {\it fa{\c c}ons}. 

\subsection{Example}

\begin{example}  \label{exal}
{\rm 

\medskip 

Let $F: \C^3 \to \C^3$ be a dominant polynomial mapping of degree 2. Assume that a 2-dimensional stratum $S$ of $S_F$  admits the two {\it fa{\c c}ons} $\kappa = (1, 2, 3)$ and $\kappa'= (1,2)[3]$. That means that all the sequences tending to infinity in the source space such that their images tend to the points of $S$ admit either the {\it fa{\c c}on} $\kappa = (1, 2, 3)$ or the {\it fa{\c c}on} $\kappa'= (1,2)[3]$. In order to describe the geometry of $S$,   we perform the following steps:


\medskip

{\bf Step 1:} Determine the pertinent variable of $F$ with respect to the {\it fa{\c c}ons} $\kappa = (1, 2, 3)$ and $\kappa'= (1,2)[3]$:

\begin{enumerate}
\item[$\bullet$] With the {\it fa{\c c}on} $\kappa = (1, 2, 3)$, up to a  suiable linear change of coordinates,  the mapping $F$ admits the pertinent variables:  $x_1-x_2$, $(x_1 - x_2)x_1$, $(x_1-x_2)x_2$, $(x_1 - x_2)x_3$,  $x_1-x_3$, $(x_1 - x_3)x_1$, $(x_1-x_3)x_2$, $(x_1 - x_3)x_3$, $x_2-x_3$, $(x_2 - x_3)x_1$, $(x_2-x_3)x_2$, $(x_2 - x_3)x_3$, $x_1 - x_2x_3$, $x_2 - x_1x_3$ and $x_3 - x_1x_2$ (see definition \ref{pertinentvar}).

\item[$\bullet$] With the {\it fa{\c c}on} $\kappa' = (1,2)[3]$, as we refer to the same mapping $F$, then up to the {\it same } suiable linear  change of coordinates, the mapping $F$ admits the pertinent variables: $x_3$, $x_1x_3$, $x_2x_3$, $x_1-x_2$, $(x_1 - x_2)x_1$, $(x_1-x_2)x_2$, $(x_1 - x_2)x_3$ and  $(x_1 - x_3)x_3$. 
\end{enumerate}

Since $S$ contains both of the {\it fa{\c c}ons} $\kappa$ and $\kappa'$, then this  surface $S$ admits  $x_1-x_2$, $(x_1 - x_2)x_1$, $(x_1-x_2)x_2$, $(x_1 - x_2)x_3$ and  $(x_1 - x_3)x_3$ as pertinent variables. Let us denote by 
$$X_1 =x_1-x_2, \quad X_2 =  (x_1 - x_2)x_1, \quad X_3 = (x_1-x_2)x_2, \quad 
X_4 = (x_1 - x_2)x_3, \quad X_5 =  (x_1 - x_3)x_3.$$
We can write 
\begin{equationth} \label{equaordre3}
F = \tilde{F}(X_1, X_2, X_3, X_4, X_5).
\end{equationth}

{\bf Step 2:}
Assume that $\{ \xi_k = (x_{1,k}, x_{2,k}, x_{3, k})\}$ and $\{\xi'_k = (x'_{1,k}, x'_{2,k}, x'_{3,k})\}$ are two sequences tending to infinity with the  {\it fa{\c c}ons} $\kappa$ and $\kappa'$, respectively. 

A) Let us consider the fa{\c c}on $\kappa = (1, 2, 3)$ and its corresponding generic sequence $\{ \xi_k = (x_{1,k}, x_{2,k}, x_{3, k})\}$:

\begin{enumerate}
\item[$\bullet$] Assume that  $X_1(\xi_k) = (x_{1,k} - x_{2,k})$ tends to a non-zero complex number. Since  $\kappa = (1, 2, 3)$ then all three coordinate sequences $x_{1,k}, x_{2,k}$ and $x_{3 , k}$ tend to infinity. Hence $X_2 (\xi_k) =  (x_{1,k} - x_{2,k})x_{1,k}$, $X_3(\xi_k) = (x_{1,k}-x_{2,k})x_{2,k}$ and  
$X_4 (\xi_k) = (x_{1,k} - x_{2,k})x_{3,k}$ tend to infinity. In this case, $X_2$, $X_3$ and $X_4$ cannot be pertinent variables of $F$ anymore. Then $F$ admits only two pertinent variables  $X_1$ and $X_5$, or $F = \tilde{F}(X_1, X_5)$. We can see that the dimension of $S$ in this case  is 1,  that provides a contradiction with the fact that the dimension of $S$ is 2. Consequently, $(x_{1,k} - x_{2,k})$ tends to 0. 

\item[$\bullet$] Assume that $(x_{1,k} - x_{3,k})$  tends to a non-zero complex number. Then $X_5(\xi_k) = (x_{1,k} - x_{3,k}) x_{3,k}$ tend to infinity, hence $X_{5}$ cannot be a pertinent variable of $F$ anymore, or $F = \tilde{F}(X_1, X_2, X_3, X_4)$. We choose a {\it generic} sequence $\{\xi_k\}$ satisfying the conditions: $X_{1,k} = (x_{1,k} - x_{2,k})$ tends to zero and $(x_{1,k} - x_{3,k})$  tends to a non-zero complex number, for example, $\xi_k = (k + \alpha/k, k + \beta/k, k + \gamma)$. Then $X_2(\xi_k) =  (x_{1,k} - x_{2,k})x_{1,k}$, $X_3(\xi_k) = (x_{1,k}-x_{2,k})x_{2,k}$ and  
$X_4(\xi_k) = (x_{1,k} - x_{2,k})x_{3,k}$ tend to the same complex number $\lambda - \mu$. Combining with the fact $X_{1,k} = (x_{1,k} - x_{2,k})$ tends to zero, we conclude that the dimension of $S$ in this case is 1, that provides a contradiction with the fact that the dimension of $S$ is 2. Consequently, $(x_{1,k} - x_{3,k})$ tends to 0. 

\end{enumerate}
Then, with the  {\it fa{\c c}on} $\kappa$, we have $(x_{1,k} - x_{2,k})$ and $(x_{1,k} - x_{3,k})$ tend to 0. Hence  $(x_{2,k} - x_{3,k})$ also tends to 0. Let us choose a {\it generic} sequence  $\{\xi_k\}$ satisfying these conditions, for example, the sequence $\left \{ \xi_k = (k + \alpha/k, k + \beta/k, k + \gamma/k) \right \}$. We see that  $X_{2, k} = (x_{1,k} - x_{2,k})x_{1,k}$, $X_{3, k} =(x_{1,k}-x_{2,k})x_{2,k}$ and $X_{4, k} = (x_{1,k} - x_{2,k})x_{3,k}$ tend to a same complex number  $\lambda = \alpha - \beta$.   Moreover, $X_{5, k} = (x_{1,k} - x_{3,k})x_{3,k}$  tends to  $\mu = \alpha- \gamma$.   So we have 
\begin{equationth} \label{faconkappa}
\lim_{k \to \infty} F(\xi_k) =  \tilde{F}(0,   \lambda,   \lambda,   \lambda,   \mu)
\end{equationth}

B) Let us consider now the fa{\c c}on $\kappa' = (1, 2)[ 3]$ and its corresponding generic  sequence $\{ \xi'_k = (x'_{1,k}, x'_{2,k}, x'_{3, k})\}$, we have two cases:

\begin{enumerate}

\item[$\bullet$] If $X_1(\xi'_k) = (x'_{1,k} - x'_{2,k})$ tends to 0: So  
 $X_4(\xi'_k) = (x'_{1,k} - x'_{2,k})x'_{3,k}$ tends to 0. We have $X_2(\xi'_k) = (x'_{1,k} - x'_{2,k})x'_{1,k}$ and $X_2(\xi'_k) = (x'_{1,k}-x'_{2,k})x'_{2,k}$ tend to a same complex number $\lambda'$ and $X_5(\xi'_k)= (x'_{1,k} - x'_{3,k}) x'_{3,k} $ tends to an arbitrary complex number $\mu'$. Then in this case, we have 
\begin{equationth} \label{faconkappa'1}
F(\xi'_k)  =   \tilde{F}(0,    \lambda',    \lambda',    0,    \mu').
\end{equationth}

\item[$\bullet$] If $X_1(\xi'_k) = (x'_{1,k} - x'_{2,k})$ tends to a non-zero complex number $\lambda' \in \C$: So  
 $X_2(\xi'_k)$  and $X_3(\xi_k)$ tend to infinity, thus $X_2$ and $X_3$ cannot be pertinent variables of $F$ anymore. Moreover, $X_4(\xi'_k)$ tends to 0 and $X_5(\xi'_k)$ tends to an arbitrary complex number $\mu'$. Then in this case, we have 
\begin{equationth} \label{faconkappa'2}
{\begin{matrix}
F & = & \tilde{F}(X_1,  & X_4, & X_5) \cr
F(\xi'_k)  & = &  \tilde{F}( \lambda', &   0,    &\mu').
\end{matrix}}
\end{equationth} 

\end{enumerate}

In conclusion, 
we have two cases: 
\begin{enumerate}
\item[1)] From (\ref{equaordre3}), (\ref{faconkappa}) and (\ref{faconkappa'1}), we have 
 $${\begin{matrix}
F & = & \tilde{F}(X_1, & X_2, & X_3, & X_4, & X_5) \cr
   {\displaystyle \lim_{k \to \infty}} F(\xi_k) & = & \tilde{F}(0,  & \lambda,  & \lambda,  & \lambda,  & \mu) \cr 
    {\displaystyle \lim_{k \to \infty}} F(\xi'_k) & = &  \tilde{F}(0, &   \lambda',  &  \lambda',  &  0,  &  \mu').
\end{matrix}} \eqno (*)$$

\item[2)] From (\ref{equaordre3}), (\ref{faconkappa}) and (\ref{faconkappa'2}), we have 
 $${\begin{matrix}
F & = & \tilde{F}(X_1, & X_4, & X_5) \cr
   {\displaystyle \lim_{k \to \infty}} F(\xi_k) & = & \tilde{F}(0,  & \lambda,  & \mu) \cr 
    {\displaystyle \lim_{k \to \infty}} F(\xi'_k) & = &  \tilde{F}( \lambda', &   0,  &  \mu'). 
\end{matrix}} \eqno (**)$$

\end{enumerate}

\medskip

{\bf Step 3:}  We restrict the pertinent variables in the step 2 using the three following facts: 
\begin{enumerate}
\item[$\bullet$] $\kappa$ and $\kappa'$
 are two {\it fa{\c c}ons} of the same  stratum $S$, 

\item[$\bullet$] $\dim S = 2$,
\item[$\bullet$] $F$ is dominant.
\end{enumerate}

Let us consider the  two cases (*) and (**) determined in the step 2: 
\begin{enumerate}
\item[1)] $F$ is of the form (*): 
\begin{enumerate}
\item[$\bullet$] At first, we use the fact that  $\kappa$ and $\kappa'$
 are two {\it fa{\c c}ons} of the same  stratum $S$, then 
 if $X_i$ is a pertinent variable of $F$ then both 
$X_i(\xi_k)$ and 
$X_i(\xi'_k)$ must tend to either an arbitrary complex number or zero.  

\item[$\bullet$] Since the dimension of $S$ is 2 then 
$F$ must have at least two pertinent variables $X_i$ and $X_j$ 
such that the images of the sequences 
$\xi_k$ and $\xi'_k$ by $X_i$ and $X_j$, respectively,  
 tend independently to two complex numbers. 
In this case: 
\begin{enumerate}
\item[$+$] $F$ must admit either $X_2 = (x_1 - x_2) x_1$ or $X_3 = (x_1 - x_2) x_2$ as a pertinent variable, 
\item[$+$] $F$ must admit $X_5 = (x_1 - x_3)x_3$ as a pertinent variable. 
\end{enumerate}
  
\item[$\bullet$] Since $F$ is dominant then $F$ must admit at least 3 independent pertinent variables (see lemma \ref{lemmaindependant}). 
Then in this case, $F$ must also admit $X_1 = x_1 - x_2$ as a pertinent variable.  
We see that $X_1(\xi_k)$ and $X_1(\xi'_k)$ tend to 0. 
We can say that this variable is a ``free'' pertinent variable. The role of this variable is to guarantee the fact  that $F(\C^3)$ is dense in the target space $\C^3$.

\end{enumerate}

\item[2)] $F$ is of the form (**): 
Similarly to the case 1), we can see easily that in this case $F$ can admit only $X_5$ as a pertinent variable. Then the dimension of $S$ is 1, which is a contradiction with the fact that the dimension of  $S$ is 2. 
\end{enumerate}

In conclusion, $F$ has the following form: 
$${\begin{matrix}
F & = & \tilde{F}(X_1, & X_2, & X_3,  & X_5)  \cr
 
   {\displaystyle \lim_{k \to \infty}} F(\xi_k) & = & \tilde{F}(0,  & \lambda,  & \lambda,  & \mu) \cr 
    {\displaystyle \lim_{k \to \infty}} F(\xi'_k) & = &  \tilde{F}(0, &   \lambda',  &  \lambda',  &  \mu'), 
\end{matrix}}$$

`

\medskip 

{\bf Step 4: Describe the geometry of the 2-dimensional stratum $S$}: On the one hand,  the pertinent variables $X_2$ (or $X_3$) and $X_5$ tending independently to two complex numbers have degree 2; on the other hand, the degree of $F$ is 2, then the degree of the surface $S$ with respect to the variables $\lambda$ and $\mu$ (or $\lambda'$ and $\mu'$) is 1  (notice that by degree of $S$, we mean the degree of the equation defining $S$). We conclude that $S$ is a plane.

\medskip 
\medskip

 }
\end{example}

In light of the example  \ref{exal}, we explicit now the algorithm for classifying the asymptotic sets of the non-proper dominant polynomial mappings $F=(F_1, F_2 , F_3): \C^3 \to \C^3$ of degree 2. 

\subsection{Algorithm}

\begin{algorithm} \label{algorithmordre} 
{\rm We have the five following steps:

\medskip 

{\bf Step 1}: 
\begin{enumerate}
\item[$\bullet$] Fix $l$ {\it fa{\c c}ons} $\kappa_1, \ldots, \kappa_l$ in the list (\ref{group}), where $1 \leq l \leq 19$. 
\item[$\bullet$] Determine the pertinent variables with respect to these $l$ {\it fa{\c c}ons} (knowing that they must be refered to a same mapping $F$).
\end{enumerate}

{\bf Step 2}:  
\begin{enumerate}
\item[$\bullet$] Assume that  $S$ is  a 2-dimensional stratum of $S_F$ admitting {\it only} the $l$ {\it fa{\c c}ons} $\kappa_1, \ldots, \kappa_l$ in step 1. 
\item[$\bullet$] Take  {\it generic} sequences  $\xi_k^1, \ldots, \xi_k^l$ corresponding to $\kappa_1, \ldots, \kappa_l$, respectively. 
\item[$\bullet$] Compute the limit of the images  of the sequences $\xi_k^1, \ldots, \xi_k^l$ by the  pertinent variables defined in step 1. 
\item[$\bullet$] Restrict the pertinent variables  in step 1 using the fact  $\dim S = 2$. 
\end{enumerate}

{\bf Step 3}: Restrict again the pertinent variables  in step 2 using the three following facts: 

\begin{enumerate}
\item[$\bullet$] the {\it fa{\c c}ons} $\kappa_1, \ldots, \kappa_l$ belongs to $S$: then the images of the generic sequences  $\xi_k^1, \ldots, \xi^l_k$ by the pertinent variables defined in the step 2 must tend to either an arbitrary complex number or zero,

\item[$\bullet$] $\dim S =2$: then  there are 
 at least two pertinent variables $X_i$ and $X_j$ 
such that the images of the sequences 
$\xi_k$ and $\xi'_k$ by $X_i$ and $X_j$, respectively,  
 tend independently to two complex numbers, 

\item[$\bullet$] F is dominant: then there are at least 3 independent pertinent variables (see lemma \ref{lemmaindependant}).
\end{enumerate}

{\bf Step 4}: Describe the geometry of the 2-dimensional irreductible stratum $S$ of $S_F$ in terms of  the pertinent variables obtained in the step 3.

{\bf Step 5}: Letting $l$ run from 1 to $19$ in the list (\ref{group}). 

}
\end{algorithm}

\begin{theorem}
{\rm With the algorithm \ref{algorithmordre}, we obtain the list of all possible asymptotic sets $S_F$ of   non-proper dominant polynomial mappings $F: \C^3 \to \C^3$ of degree 2. 
}
\end{theorem}
\begin{preuve} 
On the one hand, the process of the algorithm \ref{algorithmordre} is possible, since the number of the {\it fa{\c c}ons} in the list (\ref{group}) is finite (19 {\it fa{\c c}ons}).  On the other hand, by  the step 2, step 4 and step 5, we consider all the possible cases for all 2-dimensional irreductible strata of $S_F$. 
Since the dimension of $S_F$ is 2 (see theorem \ref{theoremjelonek1}),
 we get all the possible 
asymptotic sets $S_F$ of   non-proper dominant polynomial mappings $F: \C^3 \to \C^3$ of degree 2. 
\end{preuve}

\section{Results}
In this section, we use the algorithm \ref{algorithmordre} to prove the following theorem. 

\begin{theorem} \label{theothuyb} 

The asymptotic set of a non-proper dominant polynomial mapping $F: \C^3 \to \C^3$ of degree 2 is one of the five elements in the following list ${\mathcal{L}}_{S_F}^{(3,2)}$. Moreover, any element of this list can be realized as the asymptotic set of a dominant polynomial mapping $F: \C^3 \to \C^3$ of degree 2. 

\medskip 

{\bf The list ${\mathcal{L}}_{S_F}^{(3,2)}$: } \\
1) A plane. \\
2) A paraboloid.\\
3) The union of a plane $(\mathscr{P}): \quad r_1x_1 + r_2x_2 + r_3x_3 + r_4 = 0,$
and a plane of the form $(\mathscr{P'}): \quad r'_1x_1 + r'_2x_2 + r'_3x_3 + r'_4 = 0,$
where we can choose two of the three coefficients $r'_1, r'_2, r'_3$, 
then the third of them and the fourth coefficient $r'_4$ are determined. \\
4)  The union of a plane 
$(\mathscr{P}): \quad r_1x_1 + r_2x_2 + r_3x_3 + r_4 = 0$
and a paraboloid of the form 
$(\mathscr{H}): \quad r'_ix_i^2 + r'_jx_j + r'_lx_l + r'_4 = 0, \quad \{i, j, l \} = \{1, 2, 3 \},$ 
where we can choose two of the three coefficients $r'_1, r'_2, r'_3$, 
then the third of them and the fourth coefficient $r'_4$ are determined. \\
5) The union of three planes 
$$(\mathscr{P}): \quad r_1x_1 + r_2x_2 + r_3x_3 + r_4 = 0,$$
$$(\mathscr{P'}): \quad r'_1x_1 + r'_2x_2 + r'_3x_3 + r'_4 = 0,$$ 
$$(\mathscr{P''}): \quad r''_1x_1 + r''_2x_2 + r''_3x_3 + r''_4 = 0,$$ 
where:  

a) for $(\mathscr{P'})$, we can choose two of the three coefficients $r'_1, r'_2, r'_3$, 
then the third of them and the fourth coefficient $r'_4$ are determined,

b) for $(\mathscr{P''}) $, we can choose two of the three coefficients $r''_1, r''_2, r''_3$, 
then the third of them and the fourth coefficient $r''_4$ are determined.

\end{theorem}

In order to prove this theorem, we need the two following lemmas.

\begin{lemma} \label{cothuy2} 
Let $F = (F_1, F_2, F_3): \C^3 \to \C^3$ be a non-proper dominant polynomial mapping of degree 2.  If $S_F$ contains a surface of degree higher than 1, then either $S_F$ is a paraboloid, or $S_F$ is the union of a paraboloid and a plane. 
\end{lemma}

\begin{preuve} 
Assume that $S_F$ contains a surface  $(\mathscr{H})$. Since $\deg F = 2$ then $1 \leq \deg (\mathscr{H}) \leq 2$. 

\medskip 

A) We prove firstly that if $S_F$ contains a surface  $(\mathscr{H})$ where $\deg (\mathscr{H}) = 2$  then  $(\mathscr{H})$ is a paraboloid. Since $\deg (\mathscr{H}) = 2$ and $\deg F = 2$ then 
 $S_F$ admits one  {\it fa{\c c}on} $\kappa$ in such a way that among the pertinent variables of $F$ with respect to the {\it fa{\c c}on} $\kappa$, there exists only one {\it free} pertinent variable. That means, one of $x_1$, $x_2$ and $x_3$ is a pertinent variable of $F$ with respect to  $\kappa$.  Without loose of generality, we assume that $S_F$ admits $\kappa = (3)[2]$ as a {\it fa{\c c}on}. 
 Assume that  $ \{\xi_k  =  (x_{1,k}, x_{2,k}, x_{3,k})\}$ is a generic sequence tending to infinity with the 
 {\it fa{\c c}on} $\kappa$ and 

i) $x_{3,k}$ tends to infinity and $x_{2,k}$ tends to 0 in 
such a way that $x_{2,k}x_{3,k}$ tends to an arbitrary complex number $\lambda$,

ii) $x_{1,k}$ tends to an arbitrary complex number $\mu$. 

We see that $x^2_{1,k}$ and $(x_{1,k} + x_{2,k}) x_{1,k}$ tend to $\mu^2$. 
 Since $\deg F = 2$ and $\deg (\mathscr{H}) = 2$, then 

i) one coordinate polynomial $F_{\eta}$, where $\eta \in \{1, 2, 3 \}$,  must contain $x_1$ as an element of degree 1,

ii) the another coordinate polynomial $F_{\eta'}$, where $\eta' \in \{1, 2, 3 \}$ and $\eta' \neq \eta$,  must contain $x_1^2$ or   $(x_1 + x_{2}) x_{1}$ as a pertinent variable. 

Assume that the equation of the surface $(\mathscr{H})$ is 
 $r_1 \alpha_1^{p_1} + r_2 \alpha_2^{p_2} + r_3 \alpha_3^{p_3} + r_4 = 0.$ 
Since   $x^2_{1,k}$ and $(x_{1,k} + x_{2,k}) x_{1,k}$ tend to the same complex number $\mu^2$, and $\deg (\mathscr{H}) = 2$, then 
there exists an unique index $i \in \{ 1, 2, 3\}$ such that 
 $r_i \neq 0$  and $p_i = 2$.  
If $r_j = 0$ or  $p_j = 0$ for all $j \neq i$, $j \in \{1, 2, 3\}$, 
then $(\mathscr{H})$ is the union of two lines. That  provides the contradiction 
with the fact that  $\deg (\mathscr{H}) = 2$. 
So, there exists $j \neq i$, $j \in \{ 1, 2, 3 \}$ such that 
$r_j \neq 0$ and  $p_j = 1.$ Consequently, the surface $(\mathscr{H})$ is a paraboloid.

\medskip 

B) We prove now that if $S_F$ contains a paraboloid then the biggest possible  $S_F$ is the union of this paraboloid and a plane. 
  Since $S_F$ contains a paraboloid then with the same choice of the {\it fa{\c c}on} $\kappa=(3)[2]$ as in   A),  
the mapping $F$ must be con\-si\-de\-red as a dominant polynomial mapping of pertinent variables $x_1, x_2, x_1x_2$ and $x_2x_3$, that means:
$$F = \widetilde{F}( x_1, x_2, x_1x_2, x_2x_3).$$
We can see easily that if $x_2$ is a pertinent variable of $F$, then $S_F$ admits only the   {\it fa{\c c}on} $\kappa$ and $S_F$ is a paraboloid. 
%
%
 Assume that $S_F$ contains another irreductible surface 
 $(\mathscr{H'})$ which is different from  $(\mathscr{H})$. 
Then $F$ must be considered as a polynomial mapping of pertinent variables $x_1, x_1x_2$ and  $x_2x_3$, that means:
\begin{equationth} \label{obseverkappa}
F = \widetilde{F}( x_1, x_1x_2, x_2x_3).
\end{equationth}
Let us consider now one {\it fa{\c c}on} $\kappa'$ of $(\mathscr{H'})$ such that $\kappa' \neq \kappa$ and let  $ \{\xi'_k = (x'_{1,k}, x'_{2,k}, x'_{3,k})\}$ be a corresponding generic sequence  of $\kappa'$.  Notice that  one coordinate of $F$ admits $x_1$ as a pertinent variable. Let us show that  $x'_{1, k}$ tends to 0. Assume that 
 $x'_{1, k}$ tends to a non-zero complex number. 
As one coordinate of $F$ admits $x_1x_2$ as a pertinent variable,  then $x'_{2, k}$ does not  tend to infinity. We have two cases:

+  If $x'_{2, k}$ tends to 0, then in order to $ \xi'_k$ tending to infinity,  $x'_{3, k}$ must tend to infinity. Hence, the {\it fa{\c c}on} $\kappa'$ is $(3)[2]$. That provides the contradiction with the fact $\kappa' \neq  \kappa$. 

+ If $x'_{2, k}$ tends to a non-zero finite complex number, since one coordinate of  $F$ admits $x_2x_3$ as factor, then $x'_{3, k}$ does not tend to infinity. That provides the contradiction with the fact that $ \xi'_k$ tends to infinity. 

Therefore,  $x'_{1, k}$ tends  to 0. We have the following possible cases:

\medskip 

1) $\kappa' = (2)[1]$: then $F$ is a polynomial mapping of the form  
$F = \widetilde{F}(x_1, x_3, x_1x_2, x_1x_3).$
 Combining with  (\ref{obseverkappa}), then  
 $F = \widetilde{F}(x_1, x_1x_2).$
Therefore, $F$  is not dominant, which provides the contradiction.

\medskip 

2) $\kappa' = (3)[1]$: then $F$ is a polynomial mapping of the form  
$F = \widetilde{F}(x_1, x_2, x_1x_2, x_1x_3).$ 
Combining with  (\ref{obseverkappa}), then 
$F = \widetilde{F}(x_1, x_1x_2).$
Therefore, $F$  is not dominant, which provides the contradiction.
 
\medskip 

3) $\kappa' = (2,3)[1]$: then $F$ is a polynomial mapping of the form  
$F = \widetilde{F}(x_1, x_1x_2, x_1x_3).$ 
Combining with  (\ref{obseverkappa}), then 
$F = \widetilde{F}(x_1, x_1x_2).$ 
Therefore, $F$  is not dominant, which provides the contradiction.

\medskip 

4) $\kappa' = (3)[1,2]$: then   
$F = \widetilde{F}(x_1, x_2, x_1x_2, x_2x_3, x_3x_1).$ 
Combining with  (\ref{obseverkappa}), we have 
$F = \widetilde{F}(x_1, x_1x_2, x_2x_3).$ 
Since $x'_{1,k}x'_{2,k}$ and $x'_{1,k}$ tend to 0, then $\dim (\mathscr{H'}) \leq 1$, that provides the contradiction. 

\medskip 

5) $\kappa' = (2)[1,3]$:  in this case, $F$ is a polynomial mapping admitting  the form 
$$F = \widetilde{F}(x_1, x_3, x_1x_2, x_2x_3, x_3x_1).$$ 
Combining with  (\ref{obseverkappa}), then 
$$F = \widetilde{F}(x_1, x_1x_2, x_2x_3).$$
 We know that $x'_{1,k}$ tends to 0. 
Assume that $x'_{1,k} x'_{2,k}$ tends to a complex number $\lambda$ and  $x'_{2,k}x'_{3,k}$ tends to  a complex number $\mu$, we have 
$$(\mathscr{H'}) = \{ (\widetilde{F}_1(0, \lambda, \mu), \widetilde{F}_2(0, \lambda, \mu), \widetilde{F}_3(0, \lambda, \mu)): \lambda, \mu \in \C\},$$
where $\widetilde{F} = (\widetilde{F}_1, \widetilde{F}_2, \widetilde{F}_3)$.  
  Since $\deg F = 2$, then the degree of $\widetilde{F}_i$ with respect to the variables $\lambda$ and $\mu$ must be 1, for all $i \in \{ 1, 2, 3\}$. Consequently, the surface $(\mathscr{H'})$ is a plane. 

\end{preuve}

\begin{lemma} \label{lemmeordren=3}

Let   $F: \C^3 \to \C^3$  be a non-proper dominant polynomial mapping of degree 2. 
Assume that $S$ is a 2-dimensional irreductible stratum of $S_F$. Then $S$ admits  at most two {\it fa{\c c}ons}. Moreover, if  $S$ admits two  {\it fa{\c c}ons}, then $S_F$ is a plane. 
\end{lemma}

\begin{preuve}
Let   $F: \C^3 \to \C^3$  be a non-proper dominant polynomial mapping of degree 2. 
Assume that $S$ is a 2-dimensional irreductible stratum of $S_F$. 

\medskip 

A) We provide firstly the list of pairs of {\it fa{\c c}ons} that $S$ can admit and we write $F$ in terms of pertinent variables in each of these cases.    
 Let us fix a pair of  {\it fa{\c c}ons} $(\kappa, \kappa')$ in the list (\ref{group}) and 
 assume that $S$ admits these two {\it fa{\c c}ons}. 
We use the steps 1, 2, 3 and 4 of the algorithm \ref{algorithmordre}. In the same way than the example \ref{exal}, we can determine the form of $F$ in terms of its pertinent variables with respect to two fixed {\it fa{\c c}ons} after using the conditions of dimension of $S$ and  the dominancy of $F$.  Letting  two {\it fa{\c c}ons}  $\kappa, \kappa'$ run in the list (\ref{group}), we get the following possiblilities:


\medskip 

\noindent 1) $(\kappa, \, \kappa') = ((1, 2, 3), (i_1, i_2)[j])$, where $\{ i_1, i_2,j \} = \{ 1,2,3\}$ and 
$$F = \tilde{F}(x_{i_1}-x_{i_2}, (x_{i_1} - x_{i_2})x_{i_1}, (x_{i_1}-x_{i_2})x_{i_2}, (x_{i_1} - x_{i_2})x_j,  (x_{i_1} - x_j)x_j).$$

\noindent 2) $(\kappa, \kappa') = \left((1, 2, 3), \, (i)[j_1, j_2]\right)$, 
where $\{ i, j_1, j_2 \} = \{ 1,2,3\}$ and 
\begin{align*}
F = \tilde{F}((x_i-x_{j_1})x_{j_1}, (x_i-x_{j_1})x_{j_2}, (x_i - x_{j_2})x_{j_1}, (x_i-x_{j_2})x_{j_2}, (x_{j_1} - x_{j_2}), (x_{j_1} - x_{j_2}) x_i).
\end{align*}

\noindent 3) $(\kappa,  \kappa') = ((1, 2) [3], \, (i)[3, j])$,  
where $\{ i, j\} = \{ 1,2\}$, and
$$F = \tilde{F}(x_3, x_jx_3, x_i x_3, (x_i-x_j)x_j).$$

\medskip 

\noindent 4) $(\kappa,  \kappa')= ((1, 2)[3], \, (3)[1,2])$ and 
$$F = \tilde{F}(x_1x_3, x_2x_3, x_1 - x_2, r_1x_1x_3 + r_2x_2x_3 + r_3(x_1-x_2)x_1 + r_4(x_1 - x_2) x_2),$$
where $r_l \in \C$, for $l = 1, \ldots, 4,$ 
 such that $(r_1 \neq 0, r_2 \neq 0, (r_3, r_4) \neq (0, 0))$, or $((r_1, r_2) \neq (0, 0), (r_3, r_4) \neq (0, 0)).$

\medskip 

\noindent 5) $(\kappa, \, \kappa') = ((1, 3) [2], \,(i)[2, j])$,  where $\{i, j\}  =\{1, 3\}$, and 
$$F = \tilde{F}(x_2, x_2x_j, x_i x_2, (x_i-x_j)x_j).$$

\noindent 6) $(\kappa, \kappa')  = ((1, 3)[2], \, (2)[1,3])$ and 
$$F = \tilde{F}(x_1x_2, x_2x_3, x_1 - x_3, r_1x_1x_2 + r_2x_2x_3 + r_3(x_1-x_3)x_1 + r_4(x_1 - x_3) x_3),$$
where $r_l \in \C$, for $l = 1, \ldots, 4,$ 
such that  $(r_1 \neq 0, r_2 \neq 0)$, or $(r_1, r_2) \neq (0, 0)$ and $(r_3, r_4) \neq (0, 0)$.

\medskip 

\noindent 7) $(\kappa,  \kappa')= ((2, 3) [1], \, (i)[1, j])$, where $\{ i, j\} = \{ 2, 3\}$, and 
$$F = \tilde{F}(x_1, x_1x_i, x_1 x_j, (x_i-x_j)x_j).$$

\noindent 8) $(\kappa,  \kappa') = ((2, 3)[1], \, (1)[2,3])$ and 
$$F = \tilde{F}(x_1x_3, x_1x_2, x_3 - x_2, r_1x_1x_3 + r_2x_1x_2 + r_3(x_3-x_2)x_3 + r_4(x_3 - x_2) x_2),$$
where $r_l \in \C$, for $l = 1, \ldots, 4,$ such that  $(r_1 \neq 0, r_2 \neq 0, (r_3, r_4) \neq (0, 0))$ or $((r_1, r_2) \neq (0, 0), (r_3, r_4) \neq (0, 0))$.

\medskip 

\noindent 9) $(\kappa, \, \kappa') = ((1)[2,3], \, (i)[1, j])$, 
where $\{ i, j \} = \{ 2, 3\}$, and 
$$F = \tilde{F}(x_{j}, x_1x_{i}, r_1x_{i}x_{j} + r_2x_{1}x_j),$$
where $r_1$ et $r_2$ are the non-zero complex numbers. 

\medskip

B) We prove now that $S$ admits  at most two {\it fa{\c c}ons}. We prove the result for the first case of the above possibilities: $(\kappa, \, \kappa') = ((1, 2, 3), (i_1, i_2)[j])$, where $\{ i_1, i_2,j \} = \{ 1,2,3\}$. The other cases are proved similarly. For example, assume that 
$S$ admits two {\it fa{\c c}ons} $\kappa = (1,2,3)$ and $\kappa' = (1,2)[3]$. We prove that $S$  cannot admit the third  {\it fa{\c c}on} $\kappa''$ different from $\kappa$ and $\kappa'$. 

Let $\kappa''$ be a {\it fa{\c c}on}  of $S_F$. Let us denote by $\{\xi_k'' = (x''_{1, k}, x''_{2, k}, x''_{3, k})\}$ a generic sequence corresponding to $\kappa''$.  By the example \ref{exal}, the mapping 
$F$ admits $X_1$, $X_2$ (or $X_3$), and $X_5$ as the pertinent variables, where 
$$X_1 = x_1 - x_2, \quad X_2 = (x_1 - x_2)x_1, \quad X_3 = (x_1 - x_2)x_2, \quad X_5 = (x_1 - x_3)x_3.$$ 
Without loose the generality, we can assume that $X_2$ is a pertinent variable of $F$. 
 We prove that $X_1(\xi''_k) = (x''_{1, k} - x''_{2, k})$  tends to 0. 
 Assume that $X_1(\xi''_k) = (x''_{1, k} - x''_{2, k})$ tends to a non-zero complex number. 
 Then : 
\begin{enumerate}
\item[$+ $]
If $x''_{1,k}$ tends to infinity, then $X_2(\xi''_k)$ tends to infinity, that provides a contradiction with the fact that $X_2$ is a pertinent variable of $F$. 
\item[$+ $]
If $x''_{2,k}$ tends to infinity, then $x''_{1,k}$ also tends to infinity since $(x''_{1, k} - x''_{2, k})$ tends to a non-zero complex number. That implies $X_2(\xi''_k)$ tends to infinity and this provides a contradiction with the fact that $X_2$ is a pertinent variable of $F$. 
\end{enumerate}
Hence, $x''_{1, k}$ and  $x''_{2, k}$ cannot tend to infinity. Consequently, 
 $x''_{3, k}$ must tend to infinity. Therefore, $X_5(\xi''_k) = (x''_{1, k} - x''_{3, k}) x''_{3, k}$ tends to infinity, that provides the contradiction with the fact that $X_5$ is a pertinent variable of $F$.  We conclude that  $(x''_{1, k} - x''_{2, k})$ tend to 0. 

Then we have two possibilities:

\begin{enumerate}
\item[a)] either both of  $x''_{1, k}$ and $x''_{2, k}$ tend to 0: then $X_1(\xi''_k) = (x''_{1, k} - x''_{2, k})$, $X_2(\xi''_k) = (x''_{1, k} - x''_{2, k}) x''_{1, k}$ and $X_3(\xi''_k) =(x''_{1, k} - x''_{2, k}) x''_{2, k}$ tend to 0, 
which provides the contradiction with the fact that the dimension of $S$ is 2,

\item[b)] or both of $x''_{1, k}$ and $x''_{2, k}$ tend to infinity: Since $X_5$ is a pertinent variable of $F$, then $x''_{3, k}$ tends to 0 or infinity. We conclude that the {\it fa{\c c}on} $\kappa''$ is $(1,2,3)$ or $(1,2)[3]$. 

\end{enumerate}

In conclusion, $S$ admits only the two {\it fa{\c c}ons} $\kappa = (1,2,3)$ and $\kappa' = (1,2)[3]$. 

\medskip 

C) We prove now that if there exists a 2-dimensional irreductible stratum $S$ of $S_F$ admitting   two {\it fa{\c c}ons}, then $S_F$ is a plane. Similarly to  B), we prove this fact for the first case of the possibilities in A), that means, the case of $(\kappa, \, \kappa') = ((1, 2, 3), (i_1, i_2)[j])$, where $\{ i_1, i_2,j \} = \{ 1,2,3\}$. The other cases are proved similarly. For example, assume that 
$S$ admits two {\it fa{\c c}ons} $\kappa = (1,2,3)$ and $\kappa' = (1,2)[3]$. With the same arguments than in the example \ref{exal}, the stratum $S$ is a plane. By B), the asymptotic set $S_F$ admits also  {\it only}   two {\it fa{\c c}ons} $\kappa = (1,2,3)$ and $\kappa' = (1,2)[3]$. In other words, $S_F$ and $S$ concide. We conclude that $S_F$ is a plane.

\end{preuve}

We prove now the theorem \ref{theothuyb}.

\begin{preuve}({\it The proof of  theorem \ref{theothuyb}}). 
The cases 1) and 2)  are easily achievable by the lemmas \ref{lemmeordren=3}  and  \ref{cothuy2}, respectively. 
 Let us prove the cases 3), 4) and 5). In these cases,  on the one hand, since $S_F$ contains at least two irreductible surfaces, then $S_F$ admits at least two {\it fa{\c c}ons}; on the other hand, by the lemma \ref{lemmeordren=3}, 
each irreductible surface of $S_F$ admits only one {\it fa{\c c}on.} Assume that  $\kappa$,  $\kappa'$ are two  different {\it fa{\c c}ons} of $S_F$ and $ \{\xi_k\} $, $\{\xi'_k \}$ are  two corresponding  generic sequences, respectively.  
  We use the algorithm  \ref{algorithmordre} and in the same way than the proofs of the lemmas     \ref{cothuy2} and \ref{lemmeordren=3}, we can see easily that the pairs of  {\it fa{\c c}ons} $(\kappa, \kappa')$ must belong to only the following pairs of groups:  (I, IV), (I, V), (I, VI), (II, VI), (IV, V), (IV, VI), (V, VI) and (VI, VI) in the list   (\ref{group}).

\medskip

i) If $\kappa$ belongs to the group I and $\kappa'$ belongs to the group IV, for example $\kappa = (1, 2, 3)$ and $\kappa'= (1,2)[3]$. From the example  \ref{exal},  $F$  is a dominant polynomial mapping which can be written in terms of pertinent variables: 
$${\begin{matrix}
F &= &\tilde{F}(x_1-x_2, & (x_1 - x_2)x_1, &(x_1-x_2)x_2, &(x_1 - x_2)x_3,  &(x_1 - x_3)x_3) \cr
   {\displaystyle \lim_{k \to \infty}} F(\xi_k) & = & \tilde{F}(0, & \lambda, & \lambda, & \lambda, & \mu), \cr 
    {\displaystyle \lim_{k \to \infty}} F(\xi'_k) & = &  \tilde{F}(0, & \lambda', & \lambda', & 0, & \mu'),
\end{matrix}}$$ 
where $\lambda, \mu, \lambda', \mu' \in \C$ (see (\ref{equaordre3}), (\ref{faconkappa}) and (\ref{faconkappa'1})).  
We see that, with the sequence $\{\xi_k\}$, 
 the pertinent variables tending to an arbitrary complex numbers have the degree 2, then 
the {\it fa{\c c}on} $\kappa$ provides a plane, since the degree of $F$ is 2.
 In the same way, the {\it fa{\c c}on} $\kappa'$ provides a plane. Furthermore, it is easy to check that these two planes must have the form of the case 3) of the theorem and    $S_F$ is the  union of these two planes.

\medskip

ii) If $\kappa$ belongs to the group I  and $\kappa'$ belongs to the group V, for example  $\kappa = (1, 2, 3)$ and $\kappa'= (1)[2]$. then, on the one hand, $F$ is a dominant polynomial mapping which can be written in terms of pertinent variables:   
 \begin{equation*} \label{etyna1}
 F = \tilde{F}(x_2-x_3, (x_2 - x_3)x_2, (x_2-x_3)x_3, (x_1 - x_2)x_2).
 \end{equation*} 
On the other hand, with the same arguments than the example \ref{exal}, and for suitable generic sequences $\{\xi_k\}$ and $\{\xi'_k\}$, we obtain: 

$${\begin{matrix}
   {\displaystyle \lim_{k \to \infty}} F(\xi_k) & = & \tilde{F}(0, & \lambda, & \lambda, & \mu), \cr 
    {\displaystyle \lim_{k \to \infty}} F(\xi'_k) & = &  \tilde{F}(\lambda', & 0,  & \lambda'^2,  & \mu'),
\end{matrix}}$$
where $\lambda, \mu, \lambda', \mu' \in \C.$  
With the same arguments than in the case i), we have:

a) either the  {\it fa{\c c}ons} $\kappa$ and $\kappa'$ provide two planes of the form of the case  3) of our theorem,

b) or the {\it fa{\c c}on} $\kappa$ provides the plane $(\mathscr{P})$ and, by the lemma  \ref{cothuy2},  the {\it fa{\c c}on} $\kappa'$ provides  the paraboloid $(\mathscr{H})$ of the form of the case 4) of our theorem.  

\noindent By an easy calculation, we see that if $S_F$ admits another  {\it fa{\c c}on}  $\kappa''$ which is different from  the {\it fa{\c c}ons} $\kappa$ and  $\kappa'$, then this {\it fa{\c c}on} provides a 1-dimensional stratum contained in $(\mathscr{P})$ or contained in $(\mathscr{H}).$

\medskip

iii) Proceeding in the same way for the cases where  $(\kappa, \kappa')$ is a pair of  {\it fa{\c c}ons} belonging to the pairs of groups:  (I, VI), (II, VI), (IV, V), (IV, VI) and (V, VI), we obtain the case 3) or the case 4) of the theorem. 

\medskip

iv) Consider now the case where  $\kappa$ and  $\kappa'$ belong to the group VI, for example, $\kappa = (1)[2,3]$ and $\kappa' =(2)[1,3]$,  then $F$  is a dominant polynomial mapping which can be written in terms of pertinent variables: 
\begin{equation*} \label{etyna5}
F = \tilde{F}(x_3, x_1x_2, x_2x_3, x_3x_1).
\end{equation*}
With the same arguments than the example \ref{exal}, and for suitable generic sequences $\{\xi_k\}$ and $\{\xi'_k\}$, we obtain: 
$${\begin{matrix}
   {\displaystyle \lim_{k \to \infty}} F(\xi_k) & = & \tilde{F}( 0, &  \lambda,  & 0,  & \mu), \cr 
   {\displaystyle \lim_{k \to \infty}} F(\xi'_k) & = &  \tilde{F}(0,  & \lambda',  & \mu',  & 0), 
\end{matrix}}$$
 where $\lambda, \mu, \lambda', \mu' \in \C.$  
In this case, we have two possibilities: 
	
a) either $F$ admits $x_3$ as a free variable:  This case is similar to the case  i) and we have the case  3) of the theorem, 

b) or $F$ does not admit $x_3$  as a free variable, that means  
\begin{equation} \label{etyna6}
F = \tilde{F}(x_1x_2, x_2x_3, x_3x_1).
\end{equation}
In this case, $S_F$ admits one more {\it fa{\c c}on} $\kappa'' = (3)[1,2]$ such that with a corresponding suitable generic sequence $ \{\xi''_k \}$ of $\kappa''$, we have 
$${\begin{matrix}
   {\displaystyle \lim_{k \to \infty}} F(\xi_k) & = & \tilde{F}( \lambda,  & 0,  & \mu), \cr 
    {\displaystyle \lim_{k \to \infty}} F(\xi'_k) & = &  \tilde{F}(\lambda',  & \mu',  & 0), \cr 
    {\displaystyle \lim_{k \to \infty}} F(\xi''_k) & = &  \tilde{F}(0,  & \lambda'',  &\mu''),
 \end{matrix}}$$
 where $\lambda'', \mu'' \in \C.$  
 In this case, $S_F$ is the union of three planes the forms of which are as in the case 5) of the theorem. 
\end{preuve}

\section{The general case}

The algorithm \ref{algorithmordre}  can be generalized to clasify the asymptotic sets of    non-proper dominant polynomial mappings 
$F: \C^n \to \C^n$ of degree $d$ where $n \geq 3$ and $d \geq 2$ 
as the following.  
\begin{algorithm} \label{algorithmordregeneral} 
{\rm We have the six following steps:
\medskip

{\bf Step 1}: Determine the list ${\mathcal{L}}_F^{(n,d)}$ of all the possible {\it fa{\c c}ons} of $S_F$.  

\medskip

{\bf Step 2}: Fix $l$ {\it fa{\c c}ons} $\kappa_1, \ldots, \kappa_l$ in the list ${\mathcal{L}}_F^{(n,d)}$  obtained in step 1.
 Determine the pertinent variables with respect to these $l$ {\it fa{\c c}ons} (in the similar way than the definition  \ref{pertinentvar}).


\medskip

{\bf Step 3}: 

\begin{enumerate}
\item[$\bullet$] Assume that  $S$ is  a $(n-1)$-dimensional stratum of $S_F$ admiting {\it only} the  $l$ {\it fa{\c c}ons} $\kappa_1, \ldots, \kappa_l$ determined in step 1. 
\item[$\bullet$] Take  {\it generic} sequences  $\xi_k^1, \ldots, \xi_k^l$ corresponding to $\kappa_1, \ldots, \kappa_l$, respectively. 
\item[$\bullet$] Compute the limit of the images  of the sequences $\xi_k^1, \ldots, \xi_k^l$ by pertinent variables defined in step 1. 
\item[$\bullet$] Restrict the pertinent variables  defined in step 2 using the fact  $\dim S = n-1 $. 
\end{enumerate}



\medskip

{\bf Step 4}: Restrict again the pertinent variables in step 3 using the three following facts: 

\begin{enumerate}
\item[$\bullet$] all the {\it fa{\c c}ons}   $\kappa_1, \ldots, \kappa_l$ belong to $S$: then the images of the generic sequences  $\xi_k^1, \ldots, \xi^l_k$ by the pertinent variables defined in the step 2 must tend to either an arbitrary complex number or zero,

\item[$\bullet$] $\dim S =n-1$: then 
 there are at least $n-1$ pertinent variables $X_{i_1}, \ldots, X_{i_{n-1}}$ 
such that the images of the sequences 
$\{\xi_k^1\}, \ldots, \{\xi_k^l\}$ by $X_{i_1}, \ldots, X_{i_{n-1}}$, respectively,  
 tend independently to $(n-1)$ complex numbers, 
\item[$\bullet$] F is dominant: then there are at least $n$ independent pertinent variables (see lemma \ref{lemmaindependant}).

\end{enumerate}

{\bf Step 5}: Describe the geometry of the $(n-1)$-dimensional irreductible stratum $S$  in terms of the pertinent variables obtained in the step 4.

{\bf Step 6: } Let $l$ run in the list obtained in the step 1.

}  
\end{algorithm}

\begin{theorem}
{\rm With the algorithm \ref{algorithmordregeneral}, we obtain all possible asymptotic sets  of non-proper dominant polynomial mapping $F: \C^n \to \C^n$ of degree $d$.
}
\end{theorem}

\medskip 

\begin{preuve} 
In the one hand, 
by theorem \ref{theoremjelonek1}, the dimension of $S_F$ is $n-1$.  By  the step 3, step 5 and step 6, we consider all the possible cases of the all $(n-1)$-dimensional irreductible strata of $S_F$. 
Since the dimension of $S_F$ is $n-1$ (see theorem \ref{theoremjelonek1}),
 we get all the possible 
asymptotic sets $S_F$ of   non-proper dominant polynomial mappings $F: \C^n \to \C^n$ of degree $d$. 
 In the other hand, the number of all the possible {\it fa{\c c}ons} of a polynomial mapping $F: \C^n \to \C^n$ is finite, as the shown of the following lemma: 
\end{preuve}

\begin{lemma} {\cite{Thuy}} \label{pro nombre}
{\rm Let $F: \C^n \to \C^n$ be a polynomial mapping such that $S_F \neq \emptyset$. Then, the number of all possible {\it fa{\c c}ons} of $S_F$ is finite. More precisely, the maximum number of {\it fa{\c c}ons} of $S_F$ is equal to 
$$ \sum_{t=1}^n  C_t^n + {\displaystyle \sum_{t=1}^{n-1}} C_t^n + {\displaystyle \sum_{t=2}^{n-1}} A_t^{n}, $$
where 
$$C_t^n = \frac{n!}{t!(n-t)!},  \quad \quad \quad A_t^n = \frac{n!}{(n-t)!}.$$
}
\end{lemma}
\begin{preuve}  {\cite{Thuy}}
Assume that $\kappa = (i_1, \ldots , i_p)[j_1, \ldots, j_q]$ is a {\it fa{\c c}on} of $S_F$. 
We have the following cases: 

\medskip 

i) If  $\{ i_1, \ldots, i_p\} \cup \{ j_1, \ldots, j_q\} = \{ 1, \ldots , n\}$: 
we have ${\displaystyle \sum_{t=1}^n }  C_t^n$ possible {\it fa{\c c}ons}.

\medskip 

ii) If  $\{ i_1, \ldots, i_p\} \cup \{ j_1, \ldots, j_q\} \neq \{ 1, \ldots , n\}$ and $\{ j_1, \ldots, j_q\} = \emptyset$: we have ${\displaystyle \sum_{t=1}^{n-1}} C_t^n$ possible {\it fa{\c c}ons}. 

\medskip 

iii) If  $\{ i_1, \ldots, i_p\} \cup \{ j_1, \ldots, j_q\} \neq \{ 1, \ldots , n\}$ and $\{ j_1, \ldots, j_q\} \neq \emptyset$: We have ${\displaystyle \sum_{t=2}^{n-1}} A_t^{n}$ possible {\it fa{\c c}ons}. 

\medskip  

As the three cases are independent, then the maximum number of {\it fa{\c c}ons} of $S_F$ is equal to 
$${\displaystyle \sum_{t=1}^n }  C_t^n + {\displaystyle \sum_{t=1}^{n-1}} C_t^n + {\displaystyle \sum_{t=2}^{n-1}} A_t^{n}.$$ 
\end{preuve}

\begin{remark}
{\rm In the example \ref{exal} and in the proofs of the lemmas \ref{cothuy2} and  \ref{lemmeordren=3}, we use a linear change of variables to simplify the pertinent variables (so that we can work without coefficients and then we can simplify calculations). 
 This change of variables does not modify the results of the theorem \ref{theothuyb}. 
 However, in the algorithms \ref{algorithmordre} and \ref{algorithmordregeneral}, we do not need the step of linear change of variables, since the computers can work with coefficients of pertinent variables  without making the problem heavier.

}
\end{remark}

\bibliographystyle{plain}

\begin{thebibliography}{BBFGK}

\def\rond{\mathaccent"7017}
\def\omini{\raise 1ex \hbox{\sevenrm o}}

\bibitem{Jelonek1}  Z. Jelonek,    {\it The set of points at which polynomial map is not proper}, Ann. Polon. Math. 58 (1993), no. 3, 259-266.

\bibitem{Thuy}  T. B. T. Nguyen,   {\it La m\'ethode des fa{\c c}ons}, ArXiv: 1407.5239. 


\end{thebibliography}

\end{document}